\documentclass[11pt]{amsart}
\usepackage[utf8]{inputenc}
\usepackage[T1]{fontenc}
\usepackage{amsmath,amssymb,amsthm}
\usepackage{booktabs,longtable}
\usepackage[margin=2.7cm]{geometry}
\usepackage[hidelinks]{hyperref}

\newtheorem{theorem}{Theorem}[section]
\newtheorem{lemma}[theorem]{Lemma}
\newtheorem{proposition}[theorem]{Proposition}

\newtheorem*{theoremA}{Theorem A}
\theoremstyle{definition}

\newtheorem{remark}[theorem]{Remark}

\newcommand{\F}{\mathbb F}
\newcommand{\Z}{\mathbb Z}
\newcommand{\gr}{\operatorname{gr}}
\newcommand{\op}{\mathrm{op}}

\title[Reversible nonsymmetric rings of order 256]{All unital reversible nonsymmetric rings of order 256}

\author{José María Grau Ribas}
\address{Departamento de Matemáticas, Universidad de Oviedo, Avda.\ Calvo Sotelo s/n, 33007 Oviedo, Spain}
\email{grau@uniovi.es}

\author{Antonio M. Oller-Marcén}
\address{Centro Universitario de la Defensa de Zaragoza, Ctra.\ Huesca s/n, 50090 Zaragoza, Spain}
\email{oller@unizar.es}


\begin{document}
\begin{abstract}
A ring is reversible if $ab=0$ implies $ba=0$, and symmetric if $abc=0$ implies $acb=0$. Szabo proved that the smallest order of a unital reversible nonsymmetric ring is $256$. We classify all such minimal rings: up to isomorphism there are exactly $186$ of them, $19$ of characteristic $2$ and $167$ of characteristic $4$, and there are none of characteristic $8$ or higher. All of them are local with residue field $\F_2$, the nilpotency index of the Jacobson radical is $4$ or $5$, and $120$ of them form $60$ pairs $R$, $R^{\op}$ of non-isomorphic rings. The classification is obtained by a general method of deformation of graded algebras over $\Z/2^e\Z$, whose completeness is proved, and it is checked against the known examples: $\F_2Q_8$, Szabo's example and the ten rings in the family of generalized quaternion rings studied by Grau, Oller-Marcén and Szabo. The presentations of the $186$ rings and the code are provided as supplementary material.
\end{abstract}

\maketitle

\section{Introduction}

All rings in this paper are associative with identity $1\neq0$. A ring $R$ is \emph{reversible} if $ab=0$ implies $ba=0$ for all $a,b\in R$ \cite{Cohn}, and \emph{symmetric} if $abc=0$ implies $acb=0$ for all $a,b,c\in R$ \cite{Lambek}. Every unital symmetric ring is reversible, but the converse fails: Marks \cite{Marks} observed that the group algebra $\F_2Q_8$ of the quaternion group is reversible and nonsymmetric, and later asked \cite{Marks2} whether it is minimal with this property. Szabo \cite{Sz} answered in the affirmative: a reversible nonsymmetric ring has at least $256$ elements. In the same paper he gave a second example of order $256$, not isomorphic to $\F_2Q_8$, which is neither left nor right duo. More recently, Moore and Szabo \cite{MS} determined, for every prime $p$, the least $k$ such that there is a reversible nonsymmetric ring of order $p^k$: it is $7$ if $3\mid p-1$ and $8$ otherwise; in particular, for $p=2$ the minimal order is again $256$. None of these works classifies the rings of minimal order.

In characteristic $4$, the ring $\Z_4[i,j,k]$ of Hamilton quaternions over $\Z_4$ is reversible and nonsymmetric \cite{GMO}. Later, Grau, Oller-Marcén and Szabo \cite{GOS} studied the family
\[
\mathrm{GOS}(a,b,c,d,e,f)=\Z_4\langle x,y\rangle/\bigl(x^2-a,\ y^2-b,\ yx-xy-2(c+dx+ey+fxy)\bigr),
\]
which contains the generalized Hamilton quaternions over $\Z_4$, and found exactly $10$ reversible nonsymmetric rings among its $88$ isomorphism classes.

The aim of this paper is to complete the picture.

\begin{theoremA}
Up to isomorphism there are exactly $186$ unital reversible nonsymmetric rings of order $256$. All of them are local with $R/J\cong\F_2$, where $J$ is the Jacobson radical. Moreover:
\begin{enumerate}
\item $19$ of them have characteristic $2$ and $167$ have characteristic $4$; none has characteristic $8$ or higher;
\item the sequence of dimensions $(\dim J^i/J^{i+1})_{i\ge0}$ is $(1,2,3,2)$ for $2$ of them, $(1,2,2,2,1)$ for $64$ and $(1,3,3,1)$ for $120$; in particular, the nilpotency index of $J$ is $4$ or $5$;
\item $66$ of them are isomorphic to their opposite ring, and the remaining $120$, all with $\dim J/J^2=3$, form $60$ pairs $\{R,R^{\op}\}$ of non-isomorphic rings;
\item all of them have center of order $32$ and are duo (that is, $Ra=aR$ for every $a\in R$; equivalently, every one-sided ideal is two-sided), except two (rings 7 and 20 in the tables), which have center of order $8$ and are neither left nor right duo.
\end{enumerate}
\end{theoremA}

In characteristic $2$, the rings of type $(1,2,2,2,1)$ are $\F_2Q_8$ and five more, the ring of type $(1,2,3,2)$ is Szabo's example \cite{Sz} (see also \cite[Thm.~3.11]{Sz2}), and the twelve rings of type $(1,3,3,1)$ are new. In characteristic $4$, the ten rings of \cite{GOS}, among them the Hamilton quaternions, are among the $32$ rings of type $(1,2,2,2,1)$ with $2\in J^2\setminus J^3$; the other $157$ are new. Szabo \cite{Sz} observes that the nilpotency index of a minimal example is at most $6$ and that he knows no example with index $6$; Theorem~A shows that there is none.

The proof combines structural results, mostly from \cite{Sz}, with an exhaustive computation. We organize a local ring $R$ with $R/J=\F_2$ by its associated graded ring $\gr R=\bigoplus_i J^i/J^{i+1}$, which is an $\F_2$-algebra of dimension $8$ generated in degree $1$ \emph{in every characteristic}. The element $2\in R$ lies in $J$, and its class in $\gr R$ is a central homogeneous element. We give a complete parametrization of the rings with a given associated graded ring (Proposition~\ref{prop:param}), reduce the search to linear algebra over $\Z/2^e\Z$, and classify the results up to isomorphism by an orbit method that does not rely on invariants. The classification in characteristic $2$ was obtained by the first author in an unpublished report; here it is recovered as a particular case of the method.

\section{Reduction to local rings with residue field $\F_2$}

Throughout, $R$ is a finite ring, $J=J(R)$ its Jacobson radical and, if $R$ is local, $F=R/J$ and $d_i=\dim_F J^i/J^{i+1}$. We shall use the following results of Szabo.

\begin{theorem}[Szabo \cite{Sz}]\label{thm:szabo}
\begin{enumerate}
\item \cite[Prop.~2.3]{Sz} A reversible ring is local if and only if it is indecomposable.
\item \cite[Prop.~2.4]{Sz} If $R$ is local, reversible and $J^3=0$, then $R$ is symmetric.
\item \cite[Prop.~3.2]{Sz} If $R$ is local and reversible, $d_2=1$, and either $d_1\le 2$ or $F\cong\F_2$, then $R$ is commutative.
\item \cite[Prop.~3.4]{Sz} If $R$ is local and reversible, $F\cong\F_2$ and $(d_1,d_2,d_3)=(n,2,1)$ with $J^4=0$, then $R$ is symmetric.
\item \cite[Thm.~3.1]{Sz} Every reversible ring of order less than $256$ is symmetric.
\end{enumerate}
\end{theorem}

\begin{lemma}\label{lem:uniserial}
Let $R$ be a finite local ring with $|J/J^2|=|R/J|$. Then $R$ is symmetric.
\end{lemma}

\begin{proof}
Since $J/J^2$ is one-dimensional both as a left and as a right $F$-vector space, there is $x\in J$ with $J=Rx+J^2=xR+J^2$, and since $J$ is nilpotent, $J=Rx=xR$. By induction, $J^i=Rx^i=x^iR$. If $a\in J^i\setminus J^{i+1}$, then $a=ux^i$ with $u$ a unit, and for every unit $w$ we have $x^iw=w'x^i$ with $w'$ a unit (if $w'\in J$, then $x^i\in J^{i+1}$). Hence, if $a$, $b$ and $c$ have valuations $i$, $j$ and $k$ with respect to the $J$-adic filtration, then $abc$ is a unit times $x^{i+j+k}$, and $abc=0$ if and only if $i+j+k\ge N$, where $N$ is the nilpotency index of $J$. This condition is symmetric in $a,b,c$.
\end{proof}

\begin{proposition}\label{prop:reduction}
Let $R$ be a reversible nonsymmetric ring of order $256$. Then $R$ is local, $R/J\cong\F_2$, $J^3\neq0$, $d_1\ge2$, $d_2\ge2$, and $(d_1,d_2,d_3)\neq(4,2,1)$. Consequently, the sequence $(d_1,\dots,d_D)$, where $D$ is the largest index with $J^D\neq0$, is one of the following:
\[
(2,2,1,1,1),\ (2,2,1,2),\ (2,2,2,1),\ (2,2,3),\ (2,3,1,1),\ (2,3,2),\ (2,4,1),\ (3,2,2),\ (3,2,1,1),\ (3,3,1).
\]
\end{proposition}

\begin{proof}
A finite ring is a product of indecomposable rings, and a factor of a reversible ring is reversible. If $R$ were decomposable, each factor would have order less than $256$ and would be symmetric by Theorem~\ref{thm:szabo}(5), and so would be the product. Hence $R$ is local by Theorem~\ref{thm:szabo}(1), and $F=R/J$ is a field $\F_q$ with $q\in\{2,4,16,256\}$. If $q=256$, $R$ is a field. If $q=16$, then $|J|=16$ and $J^2=0$. If $q=4$, then $|J|=4^3$, and either $J^3=0$ or $d_1=d_2=d_3=1$, in which case $R$ is symmetric by Lemma~\ref{lem:uniserial}. Therefore $q=2$ and $\sum_{i\ge1}d_i=7$. By Theorem~\ref{thm:szabo}(2), $J^3\neq0$; by Lemma~\ref{lem:uniserial}, $d_1\ge2$; by Theorem~\ref{thm:szabo}(3), $d_2\ge2$; and the case $(4,2,1)$ is excluded by Theorem~\ref{thm:szabo}(4). Since $J^{i+1}=J^iJ$, we have $d_{i+1}\le d_1d_i$. The sequences with sum $7$, $D\ge3$, $d_1\ge2$ and $d_2\ge2$ are the ten in the statement and $(4,2,1)$.
\end{proof}

\begin{remark}
We have independently checked Theorem~\ref{thm:szabo}(3) for $R/J=\F_2$ and $d_1=2$, in every characteristic, with the method of Section~\ref{sec:method}: in the only type with $d_2=1$ and $d_1=2$ that has candidate graded algebras, namely $(1,2,1,1,1,1,1)$, all reversible rings obtained are commutative. The case $(4,2,1)$ was also checked by computer by the first author in previous work, by enumerating the $66\,422$ graded algebras of that type satisfying the condition of Lemma~\ref{lem:top}; none of them yields a nonsymmetric ring.
\end{remark}

\section{The associated graded ring}

Let $R$ be local with $R/J=\F_2$ and $|R|=256$, and let $A=\gr R=\bigoplus_{i=0}^D A_i$, $A_i=J^i/J^{i+1}$. Since $2\in J$ and $2J^i\subseteq J^{i+1}$, each $A_i$ is an $\F_2$-vector space and $A$ is a (not necessarily commutative) graded $\F_2$-algebra generated by $A_1$ with $\dim A=8$, \emph{whatever the characteristic of $R$}.

\begin{lemma}\label{lem:top}
If $R$ is reversible, then for $u\in A_i$ and $v\in A_{D-i}$ we have $uv=0$ if and only if $vu=0$.
\end{lemma}

\begin{proof}
If $\tilde u\in J^i$ and $\tilde v\in J^{D-i}$ are representatives, the product $\tilde u\tilde v\in J^D$ depends only on $u$ and $v$, because $J^{D+1}=0$. The same holds for $\tilde v\tilde u$, and it suffices to apply reversibility to $\tilde u,\tilde v$.
\end{proof}

\begin{lemma}\label{lem:D3}
Assume $J^4=0$ and $R$ reversible. Then $R$ is nonsymmetric if and only if there exist $a,b,c\in A_1$ and a permutation $\sigma$ such that $abc=0$ in $A_3$ and $\sigma(a,b,c)\neq0$. In particular, whether $R$ is symmetric depends only on $\gr R$.
\end{lemma}

\begin{proof}
Suppose $abc=0$ and $acb\ne0$ in $R$. If $c$ is a unit, then $ab=0$ and, since reversible rings are semicommutative, $acb=0$. If $a$ is a unit, then $bc=0$, hence $cb=0$ and $acb=0$. If $b$ is a unit, then $(ab)c=0$ gives $c(ab)=0$, hence $ca=0$, $ac=0$ and $acb=0$. Thus $a,b,c\in J$. If one of them lies in $J^2$, all products of the three lie in $J^4=0$. Therefore $a,b,c\in J\setminus J^2$, and all products of the three in any order lie in $J^3=A_3$ and depend only on the classes in $A_1$.
\end{proof}

\begin{lemma}\label{lem:two}
Assume that $R$ does not have characteristic $2$, and let $k\ge1$ be such that $2\in J^k\setminus J^{k+1}$. The class $t$ of $2$ in $A_k$ is central in $A$ and nonzero, and $2^m\in J^{mk}$. In particular:
\begin{enumerate}
\item if $R$ has characteristic $4$ and $2k\le D$, then $t^2=0$;
\item if $R$ has characteristic $\ge8$, then $2k\le D$;
\item if $k=1$ and $d_1\le2$, then $R$ is commutative; in particular, if $R$ is noncommutative and $k=1$, then $d_1\ge3$.
\end{enumerate}
\end{lemma}

\begin{proof}
The element $2$ is central in $R$. For (1), $t^2$ is the class of $4=0$ in $A_{2k}$. For (2), $0\ne4\in J^{2k}$. For (3), let $k=1$. If $d_1=1$, then $A_1$ is spanned by the class of $2$ and $R$ is a homomorphic image of $\Z$, hence commutative. If $d_1=2$, there is $x\in J$ such that $A_1$ is spanned by the classes of $2$ and $x$. Let $S$ be the subring generated by $x$, which contains $2$. Since $A$ is generated by $A_1$, each $A_n$ is spanned by classes of monomials in $2$ and $x$, which lie in $S$; hence $S+J^{n+1}=S+J^n$ for all $n$, and since $J$ is nilpotent, $S=R$, which is commutative. (For instance, $(\Z/4\Z)[x]/(x^4)$ has $k=1$ and $d_1=2$.)
\end{proof}

\section{Deformations of a graded algebra}\label{sec:method}

Fix a type $(1,d_1,\dots,d_D)$ and a graded $\F_2$-algebra $A=\F_2\langle x_1,\dots,x_d\rangle/I$ of that type, with $d=d_1$ and $I$ homogeneous. Choose a basis $\mathcal B=\bigcup_n\mathcal B_n$ of $A$ consisting of monomials in $x_1,\dots,x_d$, with $\mathcal B_n$ of degree $n$ and $\mathcal B_0=\{1\}$, and homogeneous relations $r_1,\dots,r_s$ of degrees between $2$ and $D$, with coefficients in $\{0,1\}$, which together with all words of length $D+1$ generate $I$ (minimal modulo those words; in practice, in each degree $n$, a basis of $I_n$ modulo $I_{n-1}V+VI_{n-1}$). We read the $r_j$ and the elements of $\mathcal B$ as noncommutative polynomials with integer coefficients.

\begin{proposition}\label{prop:param}
Let $R$ be local with $R/J=\F_2$ and $\gr R\cong A$, and let $\xi_1,\dots,\xi_d\in J$ be representatives of $x_1,\dots,x_d$.
\begin{enumerate}
\item Every element of $R$ can be written uniquely as $\sum_{b\in\mathcal B}\varepsilon_b\,b(\xi)$ with $\varepsilon_b\in\{0,1\}$.
\item In particular, there are unique $f=\sum_{\deg b\ge1}\varphi_b\,b$ and $\tau_j=\sum_{\deg b>\deg r_j}c_{j,b}\,b$, with coefficients in $\{0,1\}$, such that $2=f(\xi)$ and $r_j(\xi)=\tau_j(\xi)$.
\item $R\cong\Z\langle X_1,\dots,X_d\rangle/\bigl(2-f,\ r_j-\tau_j\ (1\le j\le s),\ \text{words of length }D+1\bigr)$.
\end{enumerate}
Conversely, for any $f$ and $\tau_j$ of this form, the quotient $Q$ in (3) is a local ring with $|Q|\le256$, and $|Q|=256$ if and only if $\gr Q\cong A$ via $X_i\mapsto x_i$.
\end{proposition}

\begin{proof}
For (1), we show by descending induction on $n$ that every element of $J^n$ is uniquely of the form $\sum_{\deg b\ge n}\varepsilon_b b(\xi)$. For $n=D+1$ this is clear. Given $z\in J^n$, its class in $A_n$ is a unique $\F_2$-combination $\sum_{b\in\mathcal B_n}\varepsilon_b b$, and $z-\sum_{b\in\mathcal B_n}\varepsilon_b b(\xi)\in J^{n+1}$, to which the induction hypothesis applies. Uniqueness follows by comparing lowest-degree terms; existence and uniqueness give $|J^n|=2^{d_n+\dots+d_D}$. Part (2) follows from (1), since $2\in J$ and $r_j(\xi)\in J^{\deg r_j+1}$. For (3), let $Q$ be the quotient. The map $X_i\mapsto\xi_i$ induces an epimorphism $Q\to R$, because $R$ is generated by the $\xi_i$ (as in the proof of Lemma~\ref{lem:two}(3)). In $Q$, the ideal $\mathfrak m=(X_1,\dots,X_d)$ is nilpotent and $Q/\mathfrak m=\Z/(2-f(0))=\F_2$, so $Q$ is local with radical $\mathfrak m$; its associated graded ring is generated by the classes of the $X_i$ and satisfies the relations $r_j$, so it is a quotient of $A$, and $|Q|=|\gr Q|\le|A|=256=|R|$. This proves (3) and the converse.
\end{proof}

Thus a ring with associated graded ring $A$ is determined by a bit vector $(\varphi,c)$, and the problem reduces to deciding for which bits the quotient has order exactly $256$. For any choice of parameters, in the quotient $Q$ we have $2=f\in\mathfrak m$ and $\mathfrak m^{D+1}=0$, so $2^{D+1}=f^{D+1}=0$; hence one may work modulo $2^{D+1}$ even before knowing whether the quotient has the right order. The computation reduces to finding the order of a submodule of $(\Z/2^{D+1}\Z)^W$, where $W$ is the number of words of length $\le D$, which we do with the Howell normal form \cite{Howell,Storjohann}. Moreover:

\begin{itemize}
\item \emph{Search by levels.} For each $L\le D$, $R/J^{L+1}$ is the analogous quotient with the tails truncated at degree $L$, and it must have order $2^{1+d_1+\dots+d_L}$. Hence the bits are chosen degree by degree, discarding branches as soon as the order fails.
\item \emph{Top degree.} If $d_D=1$, the bits of degree $D$ are determined by solving an affine system over $\F_2$ (Lemma~\ref{lem:affine} below).
\item \emph{The class of $2$.} By Lemma~\ref{lem:two}, $f$ has leading part $t$ central in $A_k$; in characteristic $4$ we require $t^2=0$ if $2k\le D$, and in characteristic $\ge8$ only $2k\le D$ needs to be considered. If $k\ge2$, then $4=f^2\equiv t^2$ modulo $J^{2k+1}$; if moreover $2k+1>D$, the characteristic is $4$ if $t^2=0$ and $8$ otherwise. The only remaining case, type $(1,2,2,1,1,1)$ with $k=2$ and $t^2=0$, was computed directly: all rings obtained have characteristic $4$. If $D=3$, Lemma~\ref{lem:two}(2) leaves only $k=1$ in characteristic $\ge8$, and by Lemma~\ref{lem:two}(3) this requires $d_1=3$; but none of the candidate graded algebras of type $(1,3,3,1)$ (nor the one of type $(1,2,3,2)$) has a nonzero central element of degree $1$, as a result of the enumeration of Section~\ref{sec:graded}. Hence $k\ge2$ in the types with $D=3$, and there the characteristic is $4$.
\end{itemize}

\begin{lemma}\label{lem:affine}
Assume $d_D=1$ and that all bits of degree $<D$ have been fixed so that the truncated quotient $Q_{D-1}$ (with $J^D=0$) has order $2^{1+d_1+\dots+d_{D-1}}$. Work in the truncated free algebra $T=(\Z/2^{D+1}\Z)\langle X\rangle/(\text{words of length }D+1)$. Let $e$ be the word in $\mathcal B_D$; let $\rho_1,\dots,\rho_m$ be the relations (including $2-f$) that carry a bit $c_j$ in degree $D$, written with $c_j=0$; and let $M_{\mathrm s}\subseteq T$ be the subgroup generated by the products $u\rho v$ of all relations with $|u|+|v|\ge1$ and by the relations without a bit in degree $D$. Then the quotient with bits $c=(c_j)$ has order $256$ if and only if
\[\gamma+\textstyle\sum_j a_jc_j\equiv0\pmod2\quad\text{for all }(\gamma,a)\in\Z^{1+m}\text{ with }\gamma e+\textstyle\sum_ja_j\rho_j\in M_{\mathrm s}.\]
In particular, the set of valid bits is the solution set, possibly empty, of an affine system over $\F_2$, and it suffices to use a generating set of the submodule of such $(\gamma,a)$, which is obtained from a Howell form of the augmented matrix.
\end{lemma}

\begin{proof}
Let $M(c)=M_{\mathrm s}+\langle\rho_j+c_je\rangle$ be the subgroup defining the quotient $Q(c)$; the products $u(c_je)v$ with $|u|+|v|\ge1$ vanish in $T$, so $M(c)$ is the correct ideal. Since $e$ has length $D\ge1$, there are words $u,v$ with $uv=e$ and $|u|+|v|\ge1$, and $u(2-f)v=2e-ufv=2e$ in $T$; hence $2e\in M_{\mathrm s}$. The subgroup $M_1=M(c)+\langle e\rangle$ does not depend on $c$, and $T/M_1=Q_{D-1}$, which has the order given by hypothesis. Since $2e\in M(c)$, the index $[M_1:M(c)]$ is $1$ or $2$, and $|Q(c)|=256$ if and only if $e\notin M(c)$. Now $e\in M(c)$ if and only if there is $(\gamma,a)$ with $\gamma e+\sum_ja_j\rho_j\in M_{\mathrm s}$ and $\sum_ja_j(\rho_j+c_je)=e+(\text{an element of }M_{\mathrm s})$, that is, with $-\gamma+\sum_ja_jc_j$ odd. Since $2e\in M_{\mathrm s}$, only parity matters, which gives the condition in the statement.
\end{proof}

\section{Classification up to isomorphism}\label{sec:iso}

Two rings with the same associated graded ring may be isomorphic even if their bits differ. We separate the mathematical criterion from its algorithmic use.

\begin{lemma}\label{lem:orbit}
Let $Q$ be a ring given by a presentation as in Proposition~\ref{prop:param}, with $|Q|=256$, and let $R'$ be another ring of order $256$. If there are $\eta_1,\dots,\eta_d\in R'$ that generate $R'$ and satisfy all relations of the presentation of $Q$, then $X_i\mapsto\eta_i$ induces an isomorphism $Q\cong R'$. If $R'$ is local with $R'/J'=\F_2$, elements $\eta_i\in J'$ generate $R'$ if and only if their classes, together with the class of $2$, span $J'/J'^2$.
\end{lemma}

\begin{proof}
The map induces an epimorphism $Q\to R'$ between rings of the same order. The second statement is the argument of Lemma~\ref{lem:two}(3).
\end{proof}

\emph{Algorithmic use.} By Proposition~\ref{prop:param}, every isomorphism $R\to R'$ between two rings in the list sends the generators $\xi$ of a presentation of $R$ to a $d$-tuple $\eta$ in $J'\setminus J'^2$ with independent classes, and the bits of $R'$ read off at $\eta$ are those of $R$. Hence, for each ring $R'$ we run through all such $d$-tuples (at most $112^3$). For each algebra $A$ in the list of its type we keep the $\eta$ such that $r(\eta)\in J'^{\deg r+1}$ for every relation $r$ of $A$, that is, those inducing $\gr R'\cong A$, and for them we read off the bits of $2$ and of the tails, which are determined by Proposition~\ref{prop:param}(2). This yields at once all presentations in the list defining rings isomorphic to $R'$, each of which is isomorphic to $R'$ by Lemma~\ref{lem:orbit}. Applied to the transposed multiplication table, the same procedure identifies the opposite ring.

\section{Enumeration of the graded algebras}\label{sec:graded}

Proposition~\ref{prop:param} gives completeness of the deformations of a fixed graded algebra. It remains to see that the list of graded algebras is complete. Let $V=A_1=\F_2^d$ and $T(V)=\bigoplus_nV^{\otimes n}$. A graded algebra of type $(1,d_1,\dots,d_D)$ generated in degree $1$ is $T(V)/I$ with $I=\bigoplus_{n\ge2}I_n$ homogeneous, $I_n\subseteq V^{\otimes n}$ of codimension $d_n$ for $n\le D$, and $I_n=V^{\otimes n}$ for $n>D$.

\begin{lemma}
A family of subspaces $I_n\subseteq V^{\otimes n}$ ($n\ge2$) is a two-sided ideal if and only if $K_n:=I_{n-1}\otimes V+V\otimes I_{n-1}\subseteq I_n$ for all $n\ge3$.
\end{lemma}

\begin{proof}
Since $T(V)$ is generated by $V$, a homogeneous subspace is a two-sided ideal if and only if it is stable under multiplication by $V$ on both sides, which is the stated condition.
\end{proof}

The algorithm builds $I$ degree by degree:
\begin{enumerate}
\item \emph{Degree $2$.} $I_2$ runs over the subspaces of codimension $d_2$ of $V\otimes V$, equivalently over the surjective linear maps $V\otimes V\to\F_2^{d_2}$ in reduced row echelon form. One representative is taken from each orbit of the action of $GL(V)$ given by $g\cdot(v\otimes w)=gv\otimes gw$. No algebra is lost up to isomorphism, because a change of basis of $A_1$ is an isomorphism of graded algebras.
\item \emph{Degree $n\ge3$.} Given $I_2,\dots,I_{n-1}$, compute $K_n$ and run over \emph{all} subspaces $I_n\supseteq K_n$ of codimension $d_n$, that is, over all surjections $V^{\otimes n}/K_n\to\F_2^{d_n}$ in reduced row echelon form (none if $\dim V^{\otimes n}/K_n<d_n$). No symmetry reduction is made here, so representatives may be repeated but none is omitted.
\item Finally $I_{D+1}=V^{\otimes(D+1)}$.
\end{enumerate}
By induction on $n$ and the previous lemma, every graded algebra of the given type is isomorphic to one of those constructed. The necessary conditions of Lemmas~\ref{lem:top} and~\ref{lem:D3} are applied to each algebra, and its central homogeneous elements, the candidates for the class of $2$ (Lemma~\ref{lem:two}), are computed. As relations $r_j$ in Proposition~\ref{prop:param} we take, in each degree $n\le D$, a basis of $I_n$ modulo $K_n$; and as $\mathcal B_n$, words whose classes form a basis of $A_n$.

\section{Computation and results}

For each type in Proposition~\ref{prop:reduction} the graded algebras were enumerated as in Section~\ref{sec:graded}, and those satisfying Lemma~\ref{lem:top} and, if $D=3$, Lemma~\ref{lem:D3} were kept. For each algebra $A$ and each possible class $t$ of $2$ (or $t=0$ in characteristic $2$), all deformations of order $256$ were found. For each of them the $256\times256$ multiplication table was built, and reversibility and nonsymmetry were checked by brute force. Table~\ref{tab:search} summarizes the computation.

\begin{table}[ht]
\centering\small
\begin{tabular}{lrrrrrr}
\toprule
& \multicolumn{3}{c}{characteristic 2} & \multicolumn{3}{c}{characteristic 4}\\
\cmidrule(lr){2-4}\cmidrule(lr){5-7}
type & order 256 & reversible & nonsymmetric & order 256 & reversible & nonsymmetric\\
\midrule
$(1,2,3,2)$ & 4 & 4 & 4 & 12 & 12 & 12\\
$(1,2,2,2,1)$ & 512 & 264 & 68 & 6\,784 & 3\,832 & 960\\
$(1,3,3,1)$ & 768 & 256 & 256 & 11\,520 & 3\,840 & 3\,840\\
$(1,2,2,1,1,1)$ & 208 & 64 & 0 & 2\,608 & 960 & 0\\
$(1,2,3,1,1)$ & 256 & 64 & 0 & 2\,240 & 960 & 0\\
$(1,3,2,1,1)$ & 15\,360 & 768 & 0 & 121\,856 & 13\,056 & 0\\
\bottomrule
\end{tabular}
\caption{Number of presentations (not of classes) obtained at each stage: with quotient of order $256$ (that is, with the prescribed associated graded ring), reversible, and reversible nonsymmetric. The types $(1,2,2,1,2)$, $(1,2,2,3)$, $(1,2,4,1)$ and $(1,3,2,2)$ have no graded algebra satisfying the necessary conditions.}\label{tab:search}
\end{table}

In characteristic $\ge8$ only the types with $D\ge4$ can occur (Lemma~\ref{lem:two}(2) and the paragraph on the class of $2$ in Section~\ref{sec:method}), and characteristics $16$ and $32$ only in type $(1,3,2,1,1)$ with $k=1$. Working modulo $2^{D+1}$ one obtains reversible rings of characteristics $8$, $16$ and $32$, but none of them is nonsymmetric (Table~\ref{tab:char8}).

\begin{table}[ht]
\centering\small
\begin{tabular}{lcrr}
\toprule
type & characteristic & reversible orbits & nonsymmetric\\
\midrule
$(1,2,2,1,1,1)$ & 8 & 80 & 0\\
$(1,2,2,2,1)$ & 8 & 237 & 0\\
$(1,2,3,1,1)$ & 8 & 274 & 0\\
$(1,3,2,1,1)$ & 8 / 16 / 32 & 384 / 16 / 121 & 0\\
\bottomrule
\end{tabular}
\caption{Characteristic $\ge8$. An orbit is a class of presentations with the same graded algebra and the same class of $2$, up to change of generators.}\label{tab:char8}
\end{table}

The classification up to isomorphism (Section~\ref{sec:iso}) gives the classes in Table~\ref{tab:classes}, which proves Theorem~A. The presentations are listed in Appendix~\ref{app:tables}. We number the rings $1$--$19$ (characteristic $2$) and $20$--$186$ (characteristic $4$).

\begin{table}[ht]
\centering\small
\begin{tabular}{lcrlcccc}
\toprule
type & char. & $k$ & No. & classes & additive group & center & $R\cong R^{\op}$\\
\midrule
$(1,2,2,2,1)$ & 2 & -- & 1--6 & 6 & $\Z_2^8$ & 32 & yes\\
$(1,2,3,2)$ & 2 & -- & 7 & 1 & $\Z_2^8$ & 8 & yes\\
$(1,3,3,1)$ & 2 & -- & 8--19 & 12 & $\Z_2^8$ & 32 & no\\
$(1,2,3,2)$ & 4 & 3 & 20 & 1 & $\Z_4\times\Z_2^6$ & 8 & yes\\
$(1,2,2,2,1)$ & 4 & 2 & 21--52 & 32 & $\Z_4^4$ & 32 & yes\\
$(1,2,2,2,1)$ & 4 & 3 & 53--72 & 20 & $\Z_4^2\times\Z_2^4$ & 32 & yes\\
$(1,2,2,2,1)$ & 4 & 4 & 73--78 & 6 & $\Z_4\times\Z_2^6$ & 32 & yes\\
$(1,3,3,1)$ & 4 & 2 & 79--174 & 96 & $\Z_4^2\times\Z_2^4$ & 32 & no\\
$(1,3,3,1)$ & 4 & 3 & 175--186 & 12 & $\Z_4\times\Z_2^6$ & 32 & no\\
\bottomrule
\end{tabular}
\caption{The $186$ classes. Here $k$ is the integer with $2\in J^k\setminus J^{k+1}$.}\label{tab:classes}
\end{table}

\begin{remark}
\begin{enumerate}
\item It never happens that $2\in J\setminus J^2$. The rings of characteristic $4$ with $2\in J^D$ are $6+1+12=19$, the same number as in characteristic $2$.
\item Rings 7 and 20 are the only non-duo rings and the only ones with center of order $8$; ring 20 is the characteristic $4$ analogue of Szabo's example, with $2=x^2y$.
\item The ten reversible nonsymmetric rings of \cite{GOS} are numbers $21,23,25,26,27,28,39,40,41,42$; the Hamilton quaternions over $\Z_4$ are number $21$.
\end{enumerate}
\end{remark}

\section{Verification}

Besides the brute-force verification of each ring, we performed the following independent checks.
\begin{enumerate}
\item \emph{Known examples.} The multiplication table of $\F_2Q_8$, built directly as a group algebra, is isomorphic to ring 5, and Szabo's nonduo example \cite{Sz,Sz2} to ring 7. Of the $256$ rings of the family $\mathrm{GOS}$, exactly $32$ are reversible nonsymmetric and they fall into $10$ isomorphism classes, as stated in \cite{GOS}; all ten are in our list.
\item \emph{Ring axioms.} For each of the $186$ rings, associativity, distributivity, the identity, reversibility (all $65\,536$ pairs), nonsymmetry and the characteristic were checked on the full table. The duo property was checked by computing, for each $a\in R$, the sets $Ra$ and $aR$ and comparing them; in rings 7 and 20 there are elements $a$ with $aR\not\subseteq Ra$ and elements $b$ with $Rb\not\subseteq bR$, so they are neither left nor right duo.
\item \emph{Printed presentations.} Each presentation in the tables was converted back into a ring by a program independent of the search, checking order $256$, reversibility and nonsymmetry.
\item \emph{Two isomorphism methods.} In type $(1,2,2,2,1)$ in characteristic $4$, the orbit classification agrees with the one obtained by comparing pairs of rings through a direct search for isomorphisms ($58$ classes). In characteristic $2$ both methods give the $19$ classes.
\item \emph{Two implementations.} The search was implemented modulo $4$ with a dedicated elimination and modulo $2^e$ with the Howell form; both agree on the types where they were compared, and the affine shortcut for the top degree agrees with exhaustive search.
\end{enumerate}

\section*{Supplementary material}
The following files accompany this paper (on arXiv, as ancillary files in the directory \texttt{anc/}).
\begin{itemize}
\item \texttt{presentations.txt}: the $186$ presentations with their invariants.
\item \texttt{code/}: the Python program (enumeration of graded algebras, search for deformations, orbit classification and checks).
\end{itemize}

\section*{Declaration of competing interest}
The authors declare that they have no known competing financial interests or personal relationships that could have appeared to influence the work reported in this paper.

\section*{Data availability}
All data and code supporting the results are included as supplementary material.

\section*{Declaration of generative AI and AI-assisted technologies in the manuscript preparation process}
During the preparation of this work the authors used Claude (Anthropic) in order to write and run the computer programs, draft the text and check the references. After using this tool, the authors reviewed and edited the content as needed and take full responsibility for the content of the published article.

\appendix
\section{Tables of presentations}\label{app:tables}

In characteristic $2$, each ring is $\F_2\langle\text{generators}\rangle$ modulo the indicated relations and the words of length $D+1$. In characteristic $4$, it is $\Z_4\langle\text{generators}\rangle$ modulo the indicated relations and the words of length $D+1$; the equality $2=f$ is part of the presentation. The column \emph{op} gives the number of the opposite ring when it is not isomorphic to the ring itself.

\subsection*{Characteristic 2, two generators, type $(1,2,2,2,1)$, $J^5=0$}
{\small\begin{longtable}{rl}
\toprule No. & relations (besides $x^2y^2=0$) \\ \midrule\endhead
1 & $x^2+y^2,\ x^2+xy+yx$ \quad (graduado de $\mathbb F_2Q_8$)\\
2 & $x^2+y^2,\ x^2+xy+yx+x^3$\\
3 & $x^2+y^2,\ x^2+xy+yx+x^3+x^2y$\\
4 & $x^2+y^2,\ x^2+xy+yx+x^3y$\\
5 & $x^2+y^2,\ x^2+xy+yx+x^3+x^2y+x^3y$ \quad ($\cong\mathbb F_2Q_8$)\\
6 & $x^2+y^2+x^3,\ x^2+xy+yx+x^3+x^3y$\\
\bottomrule\end{longtable}}

\subsection*{Characteristic 2, two generators, type $(1,2,3,2)$, $J^4=0$}
{\small\begin{longtable}{rl}
\toprule No. & relations \\ \midrule\endhead
7 & $x^2+yx+y^2,\ xyx+xy^2,\ x^2y+xyx+y^2x$ \quad ($\cong$ \cite{Sz})\\
\bottomrule\end{longtable}}

\subsection*{Characteristic 2, three generators, type $(1,3,3,1)$, $J^4=0$}
{\small\begin{longtable}{rrl}
\toprule No. & op & relations \\ \midrule\endhead
8 & 14 & $x^2+y^2,\ xz+zy,\ xy+xz+zx,\ x^2+xz+yx,\ xy+yz,\ z^2$\\
9 & 15 & $x^2+y^2+x^3,\ xz+zy,\ xy+xz+zx,\ x^2+xz+yx,\ xy+yz,\ z^2$\\
10 & 16 & $x^2+y^2,\ xz+zy,\ xy+xz+zx,\ x^2+xz+yx,\ xy+yz,\ z^2+x^3$\\
11 & 17 & $x^2+y^2+x^3,\ xz+zy,\ xy+xz+zx,\ x^2+xz+yx,\ xy+yz,\ z^2+x^3$\\
12 & 18 & $x^2+zy,\ x^2+xz+z^2,\ x^2+xy+xz+zx,\ x^2+xz+yx,\ xz+y^2,\ x^2+xy+yz,\ x^3+yzy,\ x^3+zy^2$\\
13 & 19 & $x^2+zy+x^3,\ x^2+xz+z^2,\ x^2+xy+xz+zx,\ x^2+xz+yx,\ xz+y^2,\ x^2+xy+yz,\ x^3+yzy,\ x^3+zy^2$\\
14 & 8 & $x^2+y^2,\ xy+zx,\ xz+zy,\ x^2+xz+yx,\ xy+xz+yz,\ z^2$\\
15 & 9 & $x^2+y^2+x^2y,\ xy+zx,\ xz+zy,\ x^2+xz+yx,\ xy+xz+yz,\ z^2$\\
16 & 10 & $x^2+y^2,\ xy+zx,\ xz+zy,\ x^2+xz+yx,\ xy+xz+yz,\ z^2+x^2y$\\
17 & 11 & $x^2+y^2+x^2y,\ xy+zx+x^2y,\ xz+zy,\ x^2+xz+yx,\ xy+xz+yz,\ z^2+x^2y$\\
18 & 12 & $x^2+zy,\ x^2+xz+z^2,\ xy+xz+zx,\ x^2+xz+yx,\ xz+y^2,\ x^2+xy+xz+yz,\ yzy,\ zy^2$\\
19 & 13 & $x^2+zy+x^3,\ x^2+xz+z^2,\ xy+xz+zx,\ x^2+xz+yx,\ xz+y^2,\ x^2+xy+xz+yz,\ yzy,\ zy^2$\\
\bottomrule\end{longtable}}

\subsection*{Type $(1,2,3,2)$, characteristic 4}
Ring 20: $\mathbb Z_4\langle x,y\rangle/(2=x^{2}y,\ x^{2} + xy + y^{2}=0,\ x^{3}=0,\ x^{2}y + xyx + yx^{2}=0)$, with $J^4=0$.

\subsection*{Type $(1,2,2,2,1)$, characteristic 4}
Rings $21$--$78$: $\mathbb Z_4\langle x,y\rangle/(2=f,\ x^2+xy+yx=a,\ x^2+y^2=b,\ x^4)$, with $J^5=0$. For the ten previously known rings, the last column gives one choice of parameters $(a,b,c,d,e,f)$ of the family $\mathrm{GOS}$.

{\small\begin{longtable}{rcllll}
\toprule No. & $k$ & $f$ & $a$ & $b$ & GOS \\ \midrule\endhead
21 & 2 & $x^{2}$ & $0$ & $0$ & (2,2,1,0,0,1); Hamilton quaternions \\
22 & 2 & $x^{2}$ & $0$ & $x^{3}y$ &  \\
23 & 2 & $x^{2}$ & $x^{3}y$ & $0$ & (2,2,1,0,0,0) \\
24 & 2 & $x^{2} + x^{3}y$ & $0$ & $0$ &  \\
25 & 2 & $x^{2}$ & $0$ & $x^{2}y$ & (2,3,1,1,0,1) \\
26 & 2 & $x^{2}$ & $0$ & $x^{2}y + x^{3}y$ & (2,2,1,0,1,1) \\
27 & 2 & $x^{2}$ & $x^{3}y$ & $x^{2}y$ & (2,3,1,0,0,0) \\
28 & 2 & $x^{2}$ & $x^{3}y$ & $x^{2}y + x^{3}y$ & (2,2,1,0,1,0) \\
29 & 2 & $x^{2} + x^{3}y$ & $0$ & $x^{2}y$ &  \\
30 & 2 & $x^{2} + x^{3}y$ & $0$ & $x^{2}y + x^{3}y$ &  \\
31 & 2 & $x^{2} + x^{3}y$ & $x^{3}y$ & $x^{2}y$ &  \\
32 & 2 & $x^{2} + x^{3}y$ & $x^{3}y$ & $x^{2}y + x^{3}y$ &  \\
33 & 2 & $x^{2}$ & $0$ & $x^{3}$ &  \\
34 & 2 & $x^{2}$ & $0$ & $x^{3} + x^{3}y$ &  \\
35 & 2 & $x^{2}$ & $x^{3}y$ & $x^{3}$ &  \\
36 & 2 & $x^{2} + x^{3}y$ & $0$ & $x^{3}$ &  \\
37 & 2 & $x^{2} + x^{3}y$ & $0$ & $x^{3} + x^{3}y$ &  \\
38 & 2 & $x^{2} + x^{3}y$ & $x^{3}y$ & $x^{3}$ &  \\
39 & 2 & $x^{2}$ & $x^{2}y$ & $x^{2}y$ & (2,2,1,1,1,1) \\
40 & 2 & $x^{2}$ & $x^{2}y + x^{3}y$ & $x^{2}y$ & (2,3,0,0,1,0) \\
41 & 2 & $x^{2}$ & $x^{2}y + x^{3}y$ & $x^{2}y + x^{3}y$ & (2,2,1,1,1,0) \\
42 & 2 & $x^{2} + x^{3}y$ & $x^{2}y + x^{3}y$ & $x^{2}y + x^{3}y$ & (3,3,0,1,1,1) \\
43 & 2 & $x^{2}$ & $x^{2}y$ & $x^{3} + x^{2}y$ &  \\
44 & 2 & $x^{2}$ & $x^{2}y + x^{3}y$ & $x^{3} + x^{2}y$ &  \\
45 & 2 & $x^{2}$ & $x^{2}y + x^{3}y$ & $x^{3} + x^{2}y + x^{3}y$ &  \\
46 & 2 & $x^{2} + x^{3}y$ & $x^{2}y$ & $x^{3} + x^{2}y$ &  \\
47 & 2 & $x^{2} + x^{3}y$ & $x^{2}y + x^{3}y$ & $x^{3} + x^{2}y$ &  \\
48 & 2 & $x^{2} + x^{3}y$ & $x^{2}y + x^{3}y$ & $x^{3} + x^{2}y + x^{3}y$ &  \\
49 & 2 & $x^{2} + x^{2}y$ & $0$ & $x^{3} + x^{2}y$ &  \\
50 & 2 & $x^{2} + x^{2}y$ & $x^{3}y$ & $x^{3} + x^{2}y$ &  \\
51 & 2 & $x^{2} + x^{2}y$ & $x^{3}y$ & $x^{3} + x^{2}y + x^{3}y$ &  \\
52 & 2 & $x^{2} + x^{2}y + x^{3}y$ & $x^{3}y$ & $x^{3} + x^{2}y + x^{3}y$ &  \\
53 & 3 & $x^{3}$ & $0$ & $0$ &  \\
54 & 3 & $x^{3}$ & $0$ & $x^{3}y$ &  \\
55 & 3 & $x^{3}$ & $x^{3}y$ & $0$ &  \\
56 & 3 & $x^{3} + x^{3}y$ & $0$ & $0$ &  \\
57 & 3 & $x^{3} + x^{3}y$ & $0$ & $x^{3}y$ &  \\
58 & 3 & $x^{3} + x^{3}y$ & $x^{3}y$ & $0$ &  \\
59 & 3 & $x^{3}$ & $0$ & $x^{2}y$ &  \\
60 & 3 & $x^{3}$ & $0$ & $x^{2}y + x^{3}y$ &  \\
61 & 3 & $x^{3}$ & $x^{3}y$ & $x^{2}y$ &  \\
62 & 3 & $x^{3}$ & $x^{3}y$ & $x^{2}y + x^{3}y$ &  \\
63 & 3 & $x^{3} + x^{3}y$ & $0$ & $x^{2}y$ &  \\
64 & 3 & $x^{3} + x^{3}y$ & $0$ & $x^{2}y + x^{3}y$ &  \\
65 & 3 & $x^{3} + x^{3}y$ & $x^{3}y$ & $x^{2}y$ &  \\
66 & 3 & $x^{3} + x^{3}y$ & $x^{3}y$ & $x^{2}y + x^{3}y$ &  \\
67 & 3 & $x^{3}$ & $0$ & $x^{3}$ &  \\
68 & 3 & $x^{3}$ & $0$ & $x^{3} + x^{3}y$ &  \\
69 & 3 & $x^{3}$ & $x^{3}y$ & $x^{3}$ &  \\
70 & 3 & $x^{3} + x^{3}y$ & $0$ & $x^{3}$ &  \\
71 & 3 & $x^{3} + x^{3}y$ & $0$ & $x^{3} + x^{3}y$ &  \\
72 & 3 & $x^{3} + x^{3}y$ & $x^{3}y$ & $x^{3}$ &  \\
73 & 4 & $x^{3}y$ & $0$ & $0$ &  \\
74 & 4 & $x^{3}y$ & $0$ & $x^{3}y$ &  \\
75 & 4 & $x^{3}y$ & $0$ & $x^{2}y$ &  \\
76 & 4 & $x^{3}y$ & $x^{3}y$ & $x^{2}y$ &  \\
77 & 4 & $x^{3}y$ & $0$ & $x^{3}$ &  \\
78 & 4 & $x^{3}y$ & $0$ & $x^{3} + x^{3}y$ &  \\
\bottomrule\end{longtable}}

\subsection*{Type $(1,3,3,1)$, characteristic 4, graded algebra A}
$\mathbb Z_4\langle x,y,z\rangle/(2=f,\ xz + yx=c_{1},\ x^{2} + xz + y^{2}=c_{2},\ xy + yz=0,\ xy + xz + zx=c_{4},\ x^{2} + xy + zy=c_{5},\ x^{2} + z^{2}=c_{6})$, with $J^4=0$.

{\small\begin{longtable}{rrlllllll}
\toprule No. & op & $k$ & $f$ & $c_{1}$ & $c_{2}$ & $c_{4}$ & $c_{5}$ & $c_{6}$ \\ \midrule\endhead
79 & 80 & 2 & $x^{2}$ & $0$ & $0$ & $0$ & $0$ & $0$ \\
81 & 82 & 2 & $x^{2}$ & $0$ & $0$ & $0$ & $0$ & $x^{3}$ \\
83 & 84 & 2 & $x^{2}$ & $0$ & $0$ & $0$ & $x^{3}$ & $0$ \\
85 & 86 & 2 & $x^{2}$ & $0$ & $0$ & $0$ & $x^{3}$ & $x^{3}$ \\
87 & 88 & 2 & $x^{2}$ & $0$ & $0$ & $x^{3}$ & $0$ & $0$ \\
89 & 90 & 2 & $x^{2}$ & $0$ & $0$ & $x^{3}$ & $0$ & $x^{3}$ \\
91 & 92 & 2 & $x^{2}$ & $0$ & $0$ & $x^{3}$ & $x^{3}$ & $0$ \\
93 & 94 & 2 & $x^{2}$ & $0$ & $0$ & $x^{3}$ & $x^{3}$ & $x^{3}$ \\
95 & 96 & 2 & $x^{2} + x^{3}$ & $0$ & $0$ & $0$ & $0$ & $0$ \\
97 & 98 & 2 & $x^{2} + x^{3}$ & $0$ & $0$ & $0$ & $0$ & $x^{3}$ \\
99 & 100 & 2 & $x^{2} + x^{3}$ & $0$ & $0$ & $0$ & $x^{3}$ & $0$ \\
101 & 102 & 2 & $x^{2} + x^{3}$ & $0$ & $0$ & $0$ & $x^{3}$ & $x^{3}$ \\
103 & 104 & 2 & $x^{2} + x^{3}$ & $0$ & $0$ & $x^{3}$ & $0$ & $0$ \\
105 & 106 & 2 & $x^{2} + x^{3}$ & $0$ & $0$ & $x^{3}$ & $0$ & $x^{3}$ \\
107 & 108 & 2 & $x^{2} + x^{3}$ & $0$ & $0$ & $x^{3}$ & $x^{3}$ & $0$ \\
109 & 110 & 2 & $x^{2} + x^{3}$ & $0$ & $0$ & $x^{3}$ & $x^{3}$ & $x^{3}$ \\
111 & 112 & 2 & $x^{2} + xy$ & $0$ & $0$ & $0$ & $0$ & $0$ \\
113 & 114 & 2 & $x^{2} + xy$ & $0$ & $0$ & $0$ & $0$ & $x^{3}$ \\
115 & 116 & 2 & $x^{2} + xy$ & $0$ & $0$ & $0$ & $x^{3}$ & $0$ \\
117 & 118 & 2 & $x^{2} + xy$ & $0$ & $0$ & $0$ & $x^{3}$ & $x^{3}$ \\
119 & 120 & 2 & $x^{2} + xy$ & $0$ & $0$ & $x^{3}$ & $0$ & $0$ \\
121 & 122 & 2 & $x^{2} + xy$ & $0$ & $0$ & $x^{3}$ & $0$ & $x^{3}$ \\
123 & 124 & 2 & $x^{2} + xy$ & $0$ & $0$ & $x^{3}$ & $x^{3}$ & $0$ \\
125 & 126 & 2 & $x^{2} + xy$ & $0$ & $0$ & $x^{3}$ & $x^{3}$ & $x^{3}$ \\
127 & 128 & 2 & $x^{2} + xy$ & $0$ & $x^{3}$ & $0$ & $0$ & $0$ \\
129 & 130 & 2 & $x^{2} + xy$ & $0$ & $x^{3}$ & $0$ & $0$ & $x^{3}$ \\
131 & 132 & 2 & $x^{2} + xy$ & $0$ & $x^{3}$ & $0$ & $x^{3}$ & $0$ \\
133 & 134 & 2 & $x^{2} + xy$ & $0$ & $x^{3}$ & $0$ & $x^{3}$ & $x^{3}$ \\
135 & 136 & 2 & $x^{2} + xy$ & $0$ & $x^{3}$ & $x^{3}$ & $0$ & $0$ \\
137 & 138 & 2 & $x^{2} + xy$ & $0$ & $x^{3}$ & $x^{3}$ & $0$ & $x^{3}$ \\
139 & 140 & 2 & $x^{2} + xy$ & $0$ & $x^{3}$ & $x^{3}$ & $x^{3}$ & $0$ \\
141 & 142 & 2 & $x^{2} + xy$ & $0$ & $x^{3}$ & $x^{3}$ & $x^{3}$ & $x^{3}$ \\
143 & 144 & 2 & $xy + xz$ & $0$ & $0$ & $0$ & $0$ & $0$ \\
145 & 146 & 2 & $xy + xz$ & $0$ & $0$ & $0$ & $0$ & $x^{3}$ \\
147 & 148 & 2 & $xy + xz$ & $0$ & $0$ & $0$ & $x^{3}$ & $0$ \\
149 & 150 & 2 & $xy + xz$ & $0$ & $0$ & $0$ & $x^{3}$ & $x^{3}$ \\
151 & 152 & 2 & $xy + xz$ & $x^{3}$ & $0$ & $0$ & $0$ & $0$ \\
153 & 154 & 2 & $xy + xz$ & $x^{3}$ & $0$ & $0$ & $0$ & $x^{3}$ \\
155 & 156 & 2 & $xy + xz$ & $x^{3}$ & $0$ & $x^{3}$ & $0$ & $0$ \\
157 & 158 & 2 & $xy + xz$ & $x^{3}$ & $0$ & $x^{3}$ & $0$ & $x^{3}$ \\
175 & 176 & 3 & $x^{3}$ & $0$ & $0$ & $0$ & $0$ & $0$ \\
177 & 178 & 3 & $x^{3}$ & $0$ & $0$ & $0$ & $0$ & $x^{3}$ \\
179 & 180 & 3 & $x^{3}$ & $0$ & $0$ & $0$ & $x^{3}$ & $0$ \\
181 & 182 & 3 & $x^{3}$ & $0$ & $0$ & $0$ & $x^{3}$ & $x^{3}$ \\
\bottomrule\end{longtable}}

\subsection*{Type $(1,3,3,1)$, characteristic 4, graded algebra B}
$\mathbb Z_4\langle x,y,z\rangle/(2=f,\ xz + yx=0,\ x^{2} + xy + xz + y^{2}=c_{2},\ xy + yz=0,\ xy + xz + zx=c_{4},\ x^{2} + xy + zy=c_{5},\ x^{2} + z^{2}=c_{6})$, with $J^4=0$.

{\small\begin{longtable}{rrllllll}
\toprule No. & op & $k$ & $f$ & $c_{2}$ & $c_{4}$ & $c_{5}$ & $c_{6}$ \\ \midrule\endhead
80 & 79 & 2 & $x^{2}$ & $0$ & $0$ & $0$ & $0$ \\
82 & 81 & 2 & $x^{2}$ & $x^{3}$ & $0$ & $0$ & $x^{3}$ \\
84 & 83 & 2 & $x^{2}$ & $x^{3}$ & $0$ & $x^{3}$ & $0$ \\
86 & 85 & 2 & $x^{2}$ & $0$ & $0$ & $x^{3}$ & $x^{3}$ \\
88 & 87 & 2 & $x^{2}$ & $x^{3}$ & $0$ & $0$ & $0$ \\
90 & 89 & 2 & $x^{2}$ & $0$ & $0$ & $0$ & $x^{3}$ \\
92 & 91 & 2 & $x^{2}$ & $0$ & $0$ & $x^{3}$ & $0$ \\
94 & 93 & 2 & $x^{2}$ & $x^{3}$ & $0$ & $x^{3}$ & $x^{3}$ \\
96 & 95 & 2 & $x^{2} + x^{3}$ & $0$ & $0$ & $0$ & $0$ \\
98 & 97 & 2 & $x^{2} + x^{3}$ & $x^{3}$ & $0$ & $0$ & $x^{3}$ \\
100 & 99 & 2 & $x^{2} + x^{3}$ & $x^{3}$ & $0$ & $x^{3}$ & $0$ \\
102 & 101 & 2 & $x^{2} + x^{3}$ & $0$ & $0$ & $x^{3}$ & $x^{3}$ \\
104 & 103 & 2 & $x^{2} + x^{3}$ & $x^{3}$ & $0$ & $0$ & $0$ \\
106 & 105 & 2 & $x^{2} + x^{3}$ & $0$ & $0$ & $0$ & $x^{3}$ \\
108 & 107 & 2 & $x^{2} + x^{3}$ & $0$ & $0$ & $x^{3}$ & $0$ \\
110 & 109 & 2 & $x^{2} + x^{3}$ & $x^{3}$ & $0$ & $x^{3}$ & $x^{3}$ \\
112 & 111 & 2 & $x^{2} + xy$ & $0$ & $0$ & $0$ & $0$ \\
114 & 113 & 2 & $x^{2} + xy$ & $x^{3}$ & $0$ & $0$ & $x^{3}$ \\
116 & 115 & 2 & $x^{2} + xy$ & $x^{3}$ & $0$ & $x^{3}$ & $0$ \\
118 & 117 & 2 & $x^{2} + xy$ & $0$ & $0$ & $x^{3}$ & $x^{3}$ \\
120 & 119 & 2 & $x^{2} + xy$ & $x^{3}$ & $x^{3}$ & $0$ & $0$ \\
122 & 121 & 2 & $x^{2} + xy$ & $0$ & $x^{3}$ & $0$ & $x^{3}$ \\
124 & 123 & 2 & $x^{2} + xy$ & $0$ & $x^{3}$ & $x^{3}$ & $0$ \\
126 & 125 & 2 & $x^{2} + xy$ & $x^{3}$ & $x^{3}$ & $x^{3}$ & $x^{3}$ \\
128 & 127 & 2 & $x^{2} + xy$ & $x^{3}$ & $0$ & $0$ & $0$ \\
130 & 129 & 2 & $x^{2} + xy$ & $0$ & $0$ & $0$ & $x^{3}$ \\
132 & 131 & 2 & $x^{2} + xy$ & $0$ & $0$ & $x^{3}$ & $0$ \\
134 & 133 & 2 & $x^{2} + xy$ & $x^{3}$ & $0$ & $x^{3}$ & $x^{3}$ \\
136 & 135 & 2 & $x^{2} + xy$ & $0$ & $x^{3}$ & $0$ & $0$ \\
138 & 137 & 2 & $x^{2} + xy$ & $x^{3}$ & $x^{3}$ & $0$ & $x^{3}$ \\
140 & 139 & 2 & $x^{2} + xy$ & $x^{3}$ & $x^{3}$ & $x^{3}$ & $0$ \\
142 & 141 & 2 & $x^{2} + xy$ & $0$ & $x^{3}$ & $x^{3}$ & $x^{3}$ \\
144 & 143 & 2 & $xz$ & $0$ & $0$ & $0$ & $0$ \\
146 & 145 & 2 & $xz$ & $x^{3}$ & $0$ & $0$ & $x^{3}$ \\
148 & 147 & 2 & $xz$ & $x^{3}$ & $0$ & $0$ & $0$ \\
150 & 149 & 2 & $xz$ & $0$ & $0$ & $0$ & $x^{3}$ \\
152 & 151 & 2 & $xz$ & $x^{3}$ & $x^{3}$ & $0$ & $0$ \\
154 & 153 & 2 & $xz$ & $0$ & $x^{3}$ & $0$ & $x^{3}$ \\
156 & 155 & 2 & $xz$ & $0$ & $x^{3}$ & $0$ & $0$ \\
158 & 157 & 2 & $xz$ & $x^{3}$ & $x^{3}$ & $0$ & $x^{3}$ \\
176 & 175 & 3 & $x^{3}$ & $0$ & $0$ & $0$ & $0$ \\
178 & 177 & 3 & $x^{3}$ & $x^{3}$ & $0$ & $0$ & $x^{3}$ \\
180 & 179 & 3 & $x^{3}$ & $x^{3}$ & $0$ & $0$ & $0$ \\
182 & 181 & 3 & $x^{3}$ & $0$ & $0$ & $0$ & $x^{3}$ \\
\bottomrule\end{longtable}}

\subsection*{Type $(1,3,3,1)$, characteristic 4, graded algebra C}
$\mathbb Z_4\langle x,y,z\rangle/(2=f,\ x^{2} + xy + xz + yx=0,\ x^{2} + xy + y^{2}=0,\ x^{2} + yz=0,\ x^{2} + xy + zx=c_{4},\ x^{2} + xz + zy=c_{5},\ xy + z^{2}=c_{6},\ x^{2}y=0,\ x^{3} + x^{2}z=0)$, with $J^4=0$.

{\small\begin{longtable}{rrlllll}
\toprule No. & op & $k$ & $f$ & $c_{4}$ & $c_{5}$ & $c_{6}$ \\ \midrule\endhead
159 & 160 & 2 & $x^{2}$ & $0$ & $0$ & $0$ \\
161 & 162 & 2 & $x^{2}$ & $0$ & $0$ & $x^{3}$ \\
163 & 164 & 2 & $x^{2}$ & $0$ & $x^{3}$ & $0$ \\
165 & 166 & 2 & $x^{2}$ & $x^{3}$ & $x^{3}$ & $0$ \\
167 & 168 & 2 & $x^{2} + x^{3}$ & $0$ & $0$ & $0$ \\
169 & 170 & 2 & $x^{2} + x^{3}$ & $0$ & $0$ & $x^{3}$ \\
171 & 172 & 2 & $x^{2} + x^{3}$ & $0$ & $x^{3}$ & $0$ \\
173 & 174 & 2 & $x^{2} + x^{3}$ & $x^{3}$ & $x^{3}$ & $0$ \\
183 & 184 & 3 & $x^{3}$ & $0$ & $0$ & $0$ \\
185 & 186 & 3 & $x^{3}$ & $0$ & $0$ & $x^{3}$ \\
\bottomrule\end{longtable}}

\subsection*{Type $(1,3,3,1)$, characteristic 4, graded algebra D}
$\mathbb Z_4\langle x,y,z\rangle/(2=f,\ xy + xz + yx=0,\ x^{2} + xy + y^{2}=0,\ x^{2} + yz=0,\ x^{2} + xy + zx=c_{4},\ x^{2} + xy + xz + zy=c_{5},\ xy + z^{2}=c_{6},\ x^{3} + x^{2}y=0,\ x^{2}z=0)$, with $J^4=0$.

{\small\begin{longtable}{rrlllll}
\toprule No. & op & $k$ & $f$ & $c_{4}$ & $c_{5}$ & $c_{6}$ \\ \midrule\endhead
160 & 159 & 2 & $x^{2}$ & $0$ & $0$ & $0$ \\
162 & 161 & 2 & $x^{2}$ & $0$ & $0$ & $x^{3}$ \\
164 & 163 & 2 & $x^{2}$ & $0$ & $x^{3}$ & $0$ \\
166 & 165 & 2 & $x^{2}$ & $x^{3}$ & $0$ & $0$ \\
168 & 167 & 2 & $x^{2} + x^{3}$ & $0$ & $0$ & $0$ \\
170 & 169 & 2 & $x^{2} + x^{3}$ & $0$ & $0$ & $x^{3}$ \\
172 & 171 & 2 & $x^{2} + x^{3}$ & $0$ & $x^{3}$ & $0$ \\
174 & 173 & 2 & $x^{2} + x^{3}$ & $x^{3}$ & $0$ & $0$ \\
184 & 183 & 3 & $x^{3}$ & $0$ & $0$ & $0$ \\
186 & 185 & 3 & $x^{3}$ & $0$ & $0$ & $x^{3}$ \\
\bottomrule\end{longtable}}

\end{document}